\documentclass[a4paper,11pt]{article}
\usepackage[T1]{fontenc}
\usepackage[utf8]{inputenc}
\DeclareUnicodeCharacter{03BB}{\ensuremath{\lambda}}
\usepackage[bitstream-charter]{mathdesign}

\usepackage{tikz}
\usepackage{authblk}
\usepackage{booktabs}
\usepackage{amsmath}
\usepackage{amssymb}
\usepackage{amsthm}
\usepackage{mathtools}
\theoremstyle{definition}

\usepackage{cite}
\usepackage{enumitem}
\usepackage{xcolor}
\usepackage{url}
\usepackage[all]{nowidow}
\usepackage{float}
\usepackage{multirow}
\usepackage{graphicx}
\usepackage{fancyvrb}
\usepackage{vmargin}
\setmarginsrb{2cm}{2cm}{2cm}{2cm}{0mm}{0mm}{0mm}{10mm}
\usepackage{hyperref}
\usetikzlibrary{decorations.pathmorphing,calc}
\hypersetup{colorlinks=true,allcolors=blue}
\usepackage{pst-fun}

\tikzset{%
wind turbine/.pic={
\tikzset{path/.style={fill, draw=white, ultra thick, line join=round}}
\path [path] (-.25,0) arc (180:360:.25 and .0625) -- (.0625,3) -- (-.0625,3) -- cycle;
\foreach \i in {90, 210, 330}{
\ifcase#1
\or
\path [color=blue!70!white, path, shift=(90:3), rotate=\i] 
(.5,-.1875) arc (270:90:.5 and .1875) arc (90:-90:1.5 and .1875);
\or
\path [color=red!70!white, path, shift=(90:3), rotate=\i] 
(.5,-.1875) arc (270:90:.5 and .1875) arc (90:-90:1.5 and .1875);
\or
\path [color=green!60!black, path, shift=(90:3), rotate=\i]
(.5,-.1875) arc (270:90:.5 and .1875) arc (90:-90:1.5 and .1875);
\fi
}
\path [path] (0,3) circle [radius=.25];
}}

\title{\vspace*{-1.5cm} \bfseries Mathematical modeling for sustainability: How can it promote sustainable learning in mathematics education?}
\author{Natanael Karjanto}
\affil{Department of Mathematics, University College, Natural Science Campus\par Sungkyunkwan University, Suwon~16419, Republic of Korea}

\date{\vspace*{-0.5cm} \scriptsize Updated \today}

\begin{document}
\maketitle

% Abstract
\begin{abstract}
This article reviews the current state of teaching and learning mathematical modeling in the context of sustainable development goals for education at the tertiary level. While ample research on mathematical modeling education and published textbooks on the topic are available, there is a lack of focus on mathematical modeling for sustainability. This review aims to address this gap by exploring the powerful intersection of mathematical modeling and sustainability. Mathematical modeling for sustainability connects two distinct realms: learning about the mathematics of sustainability and promoting sustainable learning in mathematics education. The former involves teaching and learning sustainability quantitatively, while the latter encompasses pedagogy that enables learners to apply quantitative knowledge and skills to everyday life and continue learning and improving mathematically beyond formal education. To demonstrate the practical application of mathematical modeling for sustainability, we discuss a specific textbook suitable for a pilot liberal arts course. We illustrate how learners can grasp mathematical concepts related to sustainability through simple yet mathematically diverse examples, which can be further developed for teaching such a course. Indeed, by filling the gap in the literature and providing practical resources, this review contributes to the advancement of mathematical modeling education in the context of sustainability. \\	

\end{abstract}

% Section 1
\section{Introduction}

Sustainable learning plays a crucial role in education as it encompasses a comprehensive approach that integrates environmental, economic, and social factors. This interdisciplinary field combines the elements of education, sociology, and environmental studies to provide learners with the tools they need to navigate a rapidly changing world. By recognizing the interconnectedness of individuals, communities, and the planet, sustainable learning strives to strike a balance that ensures long-term and continuous learning~\cite{sipos2008achieving,boeve2015the}. The goal of sustainable learning is to equip learners with the knowledge, skills, and attitudes necessary to lead a sustainable life. It goes beyond traditional education by emphasizing the importance of informed decision-making, active engagement, and responsible citizenship. Through sustainable learning, individuals become catalysts for change, contributing to the creation of a more sustainable society and better future~\cite{redman2011educating}. By integrating sustainability principles into educational practices, we can cultivate a generation of learners who are not only academically proficient, but also equipped to address the complex challenges our world faces. Through holistic education that embraces sustainable learning, individuals can be empowered to make informed choices, foster environmental stewardship, promote social equity, and drive positive change.

Sustainable learning aligns seamlessly with Sustainable Development Goal 4 (SDG4) set by the United Nations (UN). SDG4 aims to ensure inclusive and equitable quality education, fostering lifelong learning opportunities for all~\cite{un2015sgd}. Sustainable learning embodies this goal by advocating accessible education that transcends socioeconomic barriers and geographic boundaries. It emphasizes the importance of equipping learners with relevant knowledge and practical skills that can be applied in their personal and professional lives~\cite{unesco2017edu}. By embracing sustainable learning practices, we can contribute to the realization of SDG4 and its overarching vision of quality education for all. Sustainable learning recognizes that education is a powerful tool for social transformation, empowering individuals to actively participate in shaping a sustainable future. It promotes an inclusive and equitable learning environment, ensuring that every learner has equal opportunities to acquire knowledge and skills applicable in real-world contexts. Through sustainable learning, we can bridge the gap between education and the pressing challenges of our time, such as environmental degradation, social inequality, and economic instability. By integrating sustainable principles and values into educational systems, we equip learners with tools to address these complex issues and contribute to building a more sustainable and equitable world~\cite{mcgreal2017special,elfert2019lifelong}.

In line with the UN SDG4 of quality education, proficiency in mathematics holds a prominent place alongside reading skills. Mathematics education plays a pivotal role in attaining this goal, as it equips individuals with critical thinking abilities, analytical skills, and problem-solving proficiency, which are essential for success in various aspects of life. Mastering mathematics not only fosters personal growth but also contributes to sustainable economic and social development, empowering individuals to make meaningful contributions to their communities and society at large. The significance of mathematics extends beyond its fundamental nature. It provides individuals with the quantitative knowledge and critical thinking skills necessary for tackling complex issues crucial to building a sustainable future~\cite{rieckmann2017education}. Mathematical proficiency plays a vital role in addressing global warming and climate change in order to promote environmental sustainability and economic development. By nurturing a deep understanding of mathematics, learners are equipped to effectively analyze and confront these multifaceted challenges. To ensure the realization of SDG4 in quality education, it is paramount to provide high-quality mathematics education that is accessible and inclusive to all learners. This entails breaking down barriers based on socioeconomic backgrounds, gender, or ethnicity, and creating an environment that fosters mathematical excellence for all individuals. By prioritizing equitable access to mathematics education, we can enhance learners to develop the necessary skills and knowledge to actively contribute to society and advance sustainable development~\cite{cresswell2020does,li2022education}.

One of the SGD4 targets includes equal access to higher education. Although it does not explicitly mention mathematics education at the tertiary level, many subjects at the college level often encompass quantitative and analytical skills, particularly for Science, Technology, Engineering, and Mathematics (STEM) majors. These skills are essential for solid preparation in mathematics-related topics, which are integral to various fields, even in several non-STEM subjects. Some social science majors, such as economics, finance, accounting, psychology, and sociology, mathematics remains extensively utilized. These disciplines rely on mathematical models and techniques to analyze data, make predictions, and develop theories. For instance, economics employs mathematical models to understand economic behavior, predict market trends, and formulate policies. In finance and accounting, mathematics is used for analyzing financial data and making informed investment decisions. Similarly, mathematical models contribute to the understanding of human behavior and social phenomena in fields such as psychology and sociology~\cite{kline1985mathematics,brown2014mathematics,kropko2015mathematics,jacques2018mathematics,leung2020reimagining}.

The journey towards majoring in STEM disciplines often entails a sequence of mathematics courses that serve as foundational pillars. These courses range from PreCalculus to Differential Equations and encompass subjects such as Calculus, Vector Calculus, Linear Algebra, Discrete Mathematics, Engineering Mathematics, Probability and Statistics, and several other courses. However, these courses very often lean heavily toward the theoretical part rather than their practical applications. This imbalance is particularly evident for science and mathematics majors, let alone discussing sustainability in depth. Application-oriented courses are typically found within the engineering realm. As a result, many mathematics and physical science students who desire to explore broader applications of mathematics beyond their respective fields are unaware of the existence of such courses. By developing a liberal arts course centered on mathematical modeling for sustainability, our aim is twofold: to bridge this gap by introducing students to the realm of mathematical modeling and acquaint them with the mathematics underlying sustainability. Through this endeavor, we aspire to cultivate sustainable learning in mathematics education and equip students with a holistic understanding of the subject matter.

The primary objective of this article is to present a comprehensive review of research and pedagogical approaches pertaining to mathematical modeling at the undergraduate level. In addition to this, we aim to address the often-overlooked aspects of sustainability in various academic disciplines. By intertwining mathematical modeling with sustainability, we endeavor to establish a connection that facilitates learners in their exploration of sustainability and promotes the adoption of sustainable learning practices, particularly within the realm of mathematics education. This integration can be likened to the connection between the theoretical and applied aspects of mathematics through mathematical modeling. Although the existing body of literature offers numerous review articles and research papers on mathematical modeling in education, a dearth of studies focusing on the teaching and learning of mathematical modeling in the context of sustainability is evident. Therefore, this article aims to bridge this gap and contribute to the existing knowledge base by providing a much-needed review of this domain.

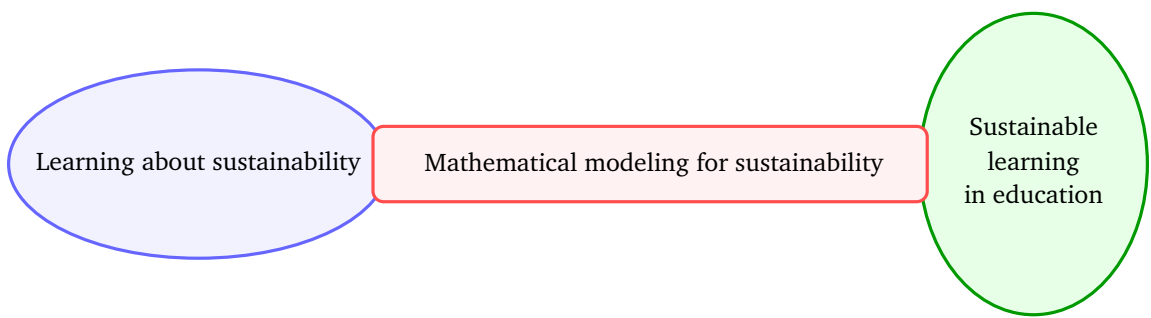
\begin{figure}[htbp]
\begin{center}
\begin{tikzpicture}
\filldraw[color=blue!60, fill=blue!5, very thick] (-6,0) ellipse (2.5 and 1.25);
\filldraw[color=green!60!black, fill=green!10, very thick](5,0) ellipse (1.5 and 2);
\filldraw[color=red!70, fill=red!5, very thick, rounded corners] (-3.7,-0.5) rectangle (3.6,0.5);
\node at (-6,0) {\small Learning about sustainability};
\node at (0,0) {\small Mathematical modeling for sustainability};
\node at (5,0.5) {\small Sustainable};
\node at (5,0) {\small learning};
\node at (5,-0.4) {\small in education};
\end{tikzpicture}
\end{center}
\caption{Illustration of the conceptual framework for this review. The blue horizontal ellipse on the left-hand side denotes the realm of ``learning about sustainability.'' This includes learning the mathematical aspects of sustainability. The green vertical ellipse on the right-hand side denotes the realm of sustainable learning in (mathematics) education, or ``learning that lasts.'' Any learner who adopts sustainable learning will continue to renew and improve themselves even after they graduate and leave the domain of formal education, particularly in mathematics and quantitative skills. The red rectangle in the middle connects these two domains, symbolizing the course on mathematical modeling for sustainability.}		\label{framework}
\end{figure}

Figure~\ref{framework} illustrates the theoretical framework of this review. This conceptual structure highlights the interconnection between learning about sustainability and sustainable learning in (mathematics) education. The left-hand side features a blue horizontal ellipse, symbolizing the realm of ``learning about sustainability.'' Here, students gain knowledge about various aspects of sustainability, including the intersection of mathematics and sustainability, as well as the relationships between the economy, culture, society, and the environment. It is our goal that through this learning, students will embrace a sustainable lifestyle throughout their lives by reusing necessary resources and minimizing wastes, actively caring for ecological health, limiting the use of Earth's natural resources, as well as ensuring quality and vitality of living environments~\cite{bhamra2007design,edwards2005the,van2013ecological,von2009factor}. 

On the right-hand side, the green vertical ellipse represents the realm of sustainable learning in (mathematics) education. This signifies learning that extends beyond the confines of formal education, emphasizing continuous renewal and improvement. Our aspiration is that students who engage in sustainable learning will consistently enhance their quantitative skills, going beyond the scope of required coursework, particularly in the mathematics and quantitative domains. They will adopt strategies in coping with challenging circumstances that require renewing and relearning, independent and collaborative learning, active learning, and lifelong learning---that is, the transferability of learning from formal to informal settings~\cite{ben2021sustainable,graham2015sustainable,hays2020sustainable,peris2015sustainable}.

Connecting these two realms is the red rectangular bridge that symbolizes the course of mathematical modeling for sustainability. Enrolling in this course offers students two primary objectives: First, it aims to provide a deeper understanding of sustainability, sustainable development, and the mathematics underlying sustainability concepts. Second, it aims to enhance learners' quantitative skills, enabling them to surpass the requirements of regular coursework. By integrating these objectives, we aspire to foster a holistic approach to sustainability and sustainable learning, empowering students to make meaningful contributions to society and the environment.

The remainder of this paper is organized as follows. Section~\ref{model} explores the world of mathematical modeling. We conducted an extensive literature review of research in education for mathematical modeling. Section~\ref{mms} continues with a liberal arts course on mathematical modeling for sustainability. We review a particularly relevant textbook and offer some additional examples that might help readers better understand the various concepts used in mathematical modeling related to sustainability. Finally, Section~\ref{discussconclude} concludes the paper.

% Section 2 % Mathematical modeling
\section{The world of mathematical modeling}	\label{model}

Mathematical modeling is the process of using mathematics to represent, analyze, and solve real-world problems by creating a mathematical representation of a scenario to make predictions or gain insights~\cite{bliss2014math}. The acquisition of principles and techniques in mathematical modeling is a powerful tool that can be utilized in various fields. Thus, teaching mathematical modeling to pupils at every level of education is essential for sustainable learning in education. As explained by Pollak (2011), the heart of mathematical modeling is problem finding before problem solving~\cite{pollak2011what}. Furthermore, Maaß et al. (2018) also argued that offering applicable mathematics education in schools will not only convince and motivate students who wonder why they should study mathematics but also equip students with transferable skills that will be useful for them along the road, such as problem-solving, critical thinking, and analytical reasoning, making it suitable for promoting sustainable learning in mathematics education~\cite{maasz2018mathematical}.

In what follows, we provide a literature review of the teaching and learning of mathematical modeling. We continue with some textbooks that can be adopted as useful resources for teaching a course on mathematical modeling at the tertiary level, while attempting to identify topics that are significantly relevant to sustainability.  

% Subsection 2.1
\subsection{Literature review}					\label{literature}

Greefrath and Vorhölter (2016) provided an overview of the German discussion on modeling and applications in schools, considering the development from the beginning of the 20th century to the present, and discussed the term ``mathematical model'' as well as different representations of the modeling process as modeling cycles~\cite{greefrath2016teaching}. Greer (1993) discovered that 13- and 14-year-old students often provide unrealistic answers to word problems because of their tendency to rely on stereotyped procedures and assumptions of direct proportionality, highlighting the need for a shift in perspective towards viewing word problems as modeling exercises, while considering the underlying assumptions and appropriateness of the model used~\cite{greer1993the}. Gainsburg (2006) reported an ethnographic study of structural engineers and revealed that modeling is central to their work. This presents challenges for K-12 education because, on the one hand, reformers in mathematics education advocated for incorporating mathematical modeling activities into K-12 curricula to reflect real-world problem solving, and on the other hand, a lack of observational descriptions of adult modeling behavior makes it difficult to assess the authenticity of classroom modeling tasks~\cite{gainsburg2006the}.

Similarly, Niss and Blum (2020), Leung et al. (2021), Saxena et al. (2016), Arseven (2015), Dunn and Marshman (2020), and Schukajlow et al. (2018) provide insights and perspectives on teaching and learning mathematical modeling in various educational contexts~\cite{niss2020the,leung2021mathematical,saxena2016teaching,arseven2015mathematical,dunn2020teaching,schukajlow2018empirical}. With an emphasis on the secondary level, Niss and Blum (2020) included resources for teachers to acquire the knowledge and competencies that will allow them to successfully include modeling in their teaching~\cite{niss2020the}. Leung et al. (2021) explored the teaching and learning of mathematical modeling and its applications, showcasing research and collaboration between colleagues from China and other parts of the world, offering new perspectives and resources for mathematical modeling education~\cite{leung2021mathematical}. Saxena et al. (2016) presented the teaching and learning of mathematical modeling from an Indian perspective, discussing the benefits and challenges while attempting to  improve the traditional way of the process~\cite{saxena2016teaching}. 

In addition to providing an overview of the theoretical basis, concepts, and inclusion of modeling activities in Turkish primary, secondary, and high school mathematics schools, Arseven (2015) explored the importance of mathematical modeling in science and mathematics education, specifically in the context of international exams such as PISA, where the report can be useful as a pathway to improving sustainable education for each participant country~\cite{arseven2015mathematical,gutierrez2022digital}. Dunn and Marshman (2020) highlighted the importance of connecting mathematics with the real world through mathematical modeling using real data, emphasized the need to understand the implications and differences between validating and estimating with data, and recommended an approach of mathematizing the context to promote students' perception of the real-world relevance of mathematics~\cite{dunn2020teaching}. Schukajlow et al. (2018) surveyed the current status of empirical studies in the teaching and learning of mathematical modeling, analyzed the development of studies focusing on cognitive aspects of the promotion of modeling, and discovered that case studies and cognitively oriented approaches tend to be more dominant than quantitative research or affect-related studies~\cite{schukajlow2018empirical}. 

Kaiser (2017) offered valuable information and perspectives for researchers and educators interested in the field of mathematical modeling in mathematics education. The article  covers a range of topics, such as the definition and characteristics of mathematical modeling, the role of teachers and students in the modeling process, challenges and strategies for teaching modeling, and the assessment of modeling skills~\cite{kaiser2017the}. In their theoretical study, Dundar et al. (2012) argued that mathematical modeling is a foundational aspect of mathematics education, involving the conversion of real-world problems into mathematical forms, and applicable across various areas of mathematics and educational levels~\cite{dundar2012mathematical}. Recently, Cevikbas et al. (2022) reviewed the current discussion on mathematical modeling competencies, urging the need for further the theoretical work, while highlighting the richness of developed empirical approaches and their implementation at various educational levels~\cite{cevikbas2022a}.

% Subsection 2.2
\subsection{Textbooks on mathematical modeling}	\label{textbook}

Many books and monographs are dedicated to mathematical modeling, ranging from an introductory level to analysis and simulations. In what follows, the list is far from exhaustive, and although they are not particularly concerned with sustainability and sustainable development, some topics and examples can still relate to, be applied to, and be useful for understanding issues in sustainability.

For example, the second part of Haberman (1988) deals with mathematical ecology and population dynamics, which sparks recent interest in sustainable populations due to low birth rates in many parts of the world~\cite{haberman1988mathematical}. When discussing population ecology, Mooney and Swift (1999) introduced the harvesting effect to the population model, raising some awareness of the sustainability of a particular population due to hunting, fishing, culling, or lumbering~\cite{mooney1999a}. In addition to discussing experimental modeling of harvesting bluefish (\emph{Pomatomus saltatrix}) and blue crabs (\emph{Callinectes sapidus}) from Chesapeake Bay in Maryland and Virgina, Giordano et al. (2014) also presented the management of renewable resources, such as the fishing industry~\cite{giordano2014a}. Heinz (2011) covered both the deterministic and stochastic aspects of global warming modeling and carbon dioxide concentration~\cite{heinz2011mathematical}. 

Although it does not contain the phrase ``mathematical modeling'' in its title, the second edition of Shiflet and Shiflet (2014) covers a wide range of mathematical modeling and its applications, where several dedicated modules are strongly related to sustainability, such as disease modeling, the carbon cycle, global warming, mercury pollution, flu pandemics, and ornamental gardens~\cite{shiflet2014introduction}. Similarly, Tung (2007) dedicated one chapter each to snowball Earth and global warming, as well as to El Niño and the Southern Oscillation (ENSO), the topics that relate to climate patterns and their consequences toward sustainability~\cite{tung2007topics}. Meerschaert (2013) presented examples related to population dynamics of two similar species that are competing between each other, that is blue whale (\emph{Balaenoptera musculus}) and fin whale (\emph{Balaenoptera physalus})~\cite{meerschaert2013mathematical}. Primarily intended for students with a working knowledge of calculus but minimal training in computer programming in the first course of modeling, Humi (2017) assists readers in mastering the processes used by scientists and engineers in modeling real-world problems, including the challenges posed by space exploration, climate change, energy sustainability, chaotic dynamical systems, and random processes~\cite{humi2017introduction}.

Other textbooks on mathematical modeling do not necessarily contain examples or applications that are directly related to sustainability. However, instructors who plan to utilize such materials can be creative and innovative in connecting mathematical topics with sustainability applications. Very often, these textbooks cover a wide range of applications, not only in physical science and engineering, but also in biological systems and social sciences, issues that have some implications for sustainability. 

For example, Meyer (1984) featured independent sections on mathematical modeling that illustrated the most important principles and a variety of applications from physical structures, biological systems, social sciences, and operational research~\cite{meyer2004concepts}. Similarly, employing a practical learning-by-doing approach, Bender (2000) offered more than 100 reality-based examples from various fields, while fostering the development of skills needed to set up and manipulate mathematical models~\cite{bender2000an}. After introducing a set of foundational tools in mathematical modeling, Dym (2004) applied these tools to a broad variety of fields that range from biology to economics, including traffic flow, mechanical systems, optimization problems, and social decision-making~\cite{dym2004principles}. Howison (2005) demonstrated that applied mathematics is much more than a series of academic calculations through a dozen applications, from the modeling of hair to piano tuning, egg incubation, and traffic flow~\cite{howison2005practical}. Using only basic knowledge of calculus and linear algebra, Velten (2009) also selected many examples from various fields, such as biology, ecology, economics, medicine, agriculture, chemical, electrical, mechanical, and process engineering~\cite{velten2009mathematical}.

The volume from Lin and Segel (1988) can be considered a classic in the pedagogy of applied mathematics, where it addressed the construction, analysis, and interpretation of mathematical models that shed light on significant problems in the physical sciences. This book is divided into three parts. The first part provides an overview of the interaction between mathematics and natural sciences. The second part illustrates some fundamental procedures in ordinary differential equations. Finally, the third part discusses the theory of continuous fields~\cite{lin1988mathematics}. Fowkes and Mahony (1994) covered a wide range of mathematics in their textbook, and explained why some techniques work well in certain situations but not in other circumstances. Using the apprenticeship approach, they also suggested when to stop analyzing mathematically and start the numerical work instead~\cite{fowkes1994an}. 

Gernshenfeld (1999) discussed common techniques that are useful for mathematical modeling, such as analytical techniques, numerical methods, and observational data analysis~\cite{gernshenfeld1999the}. Bungartz et al. (2014) not only introduced mathematical modeling but also emphasized the importance of  computer-oriented modeling and simulation as a universal methodology. In addition to addressing various model classes and their derivations, the authors also showcased the versatility of approaches that can be applied, including discrete, continuous, deterministic, and stochastic methods~\cite{bungartz2014modeling}. Eck et al. (2017) featured the use of mathematical structures as an ordering principle instead of the fields of application in their modeling book although the authors also discussed several applications in the fields of population and chemical reaction dynamics, among others~\cite{eck2017mathematical}. Garfinkel et al. (2017) specifically presented the modeling approach with examples and applications in life sciences. It equips readers with the quantitative skills needed in exploring complex feedback relationships and counterintuitive responses commonly found in natural dynamical systems~\cite{garfinkel2017modeling}. By targeting his textbook at newcomers to mathematical modeling, Banerjee (2021) offered step-by-step guidance on model formulation and covered a wide range of examples from different fields, providing an interdisciplinary overview and highlighting common themes in equilibrium, stability, bifurcations, and parameter estimation~\cite{banerjee2021mathematical}.

Many textbooks in mathematical modeling focus on continuous variables instead of discrete ones. In terms of randomness and uncertainty, they mostly focus on deterministic rather than stochastic or probabilistic aspects. Although several textbooks offer both analytical and numerical approaches, most of them focus on the former. Table~\ref{tablecharacteristics} summarizes these observations.
\begin{table}[hbtp]
\begin{center}
\begin{tabular}{@{}lccccccc@{}}
\toprule
\multicolumn{2}{c}{Textbook} 	& \multicolumn{2}{c}{Type of variables} & \multicolumn{2}{c}{Randomness and uncertainty} & \multicolumn{2}{c}{Approach}	\\ \midrule
Author(s)				& Year 	& Discrete 			& Continuous		& Deterministic			& Stochastic/ 	& Analytical 	& Numerical 			\\ 
						&		&					&					&						& probabilistic	&				&						\\	\midrule
Banerjee				& 2021	& \checkmark		& \checkmark		& \checkmark			& \checkmark	& \checkmark	& 					\\						
Bender					& 2000	& 					& \checkmark		& \checkmark			& \checkmark	& \checkmark	& \checkmark		\\
Bungartz et al. 		& 2014	& \checkmark		& \checkmark		& \checkmark			& \checkmark	& \checkmark	& \checkmark		\\
Dym		 		 		& 2004	& 					& \checkmark		& \checkmark			& 				& \checkmark	&					\\
Eck et al. 		 		& 2010	& 					& \checkmark		& \checkmark			& 				& \checkmark	&					\\
Fowkes \& Mahony 		& 1994	& 					& \checkmark		& \checkmark			& 				& \checkmark	& \checkmark		\\
Garfinkel et al.  		& 2017	& 					& \checkmark		& \checkmark			& 				& \checkmark	&					\\
Gershenfeld		 		& 2003	& 					& \checkmark		& \checkmark			& 				& \checkmark	& \checkmark		\\
Giordano et al. 		& 2014	& \checkmark		& \checkmark		& \checkmark			& \checkmark	& \checkmark	&					\\
Haberman 		 		& 1988	& 					& \checkmark		& \checkmark			& 				& \checkmark	&					\\
Heinz	 		 		& 2011	& 					& \checkmark		& \checkmark			& \checkmark	& \checkmark	&					\\
Howison 		 		& 2010	& 					& \checkmark		& \checkmark			& 				& \checkmark	&					\\
Humi	 		 		& 2017	& 					& \checkmark		& \checkmark			& 				& \checkmark	&					\\
Lin \& Segel	 		& 1988	& \checkmark		& \checkmark		& \checkmark			& 				& \checkmark	&					\\
Meerschaert		 		& 2013	& \checkmark		& \checkmark		& \checkmark			& \checkmark	& \checkmark	&					\\
Meyer			 		& 2004	& \checkmark		& \checkmark		& \checkmark			& \checkmark	& \checkmark	&					\\
Mooney \& Swift 		& 1999	& \checkmark		& \checkmark		& \checkmark			& \checkmark	& \checkmark	&					\\
Shiflet \& Shiflet 		& 2014	& 					& \checkmark		& \checkmark			& 				& \checkmark	& \checkmark		\\
Tung	 		 		& 2007	& 					& \checkmark		& \checkmark			& 				& \checkmark	&					\\
Velten	 		 		& 2009	& 					& \checkmark		& \checkmark			& 				& \checkmark	& \checkmark		\\
\bottomrule
\end{tabular}
\end{center}
\caption{Summary of textbook content based on three main characteristics. First, whether the type of variables involved in the modeling is discrete or continuous. Second, whether the factors of randomness and uncertainty were considered, that is, deterministic vs. stochastic or probabilistic. Third, whether the approach to solving modeling problems is analytical or numerical/computational.}		\label{tablecharacteristics}
\end{table}

\section{Mathematical modeling for sustainability}	\label{mms}

To the best of our knowledge, there are currently two textbooks dealing with mathematical modeling for sustainability: Roe et al. (2018)~\cite{roe2018mathematics} and Hersh (2006)~\cite{hersh2006mathematical}. Other books have some potential to be used for teaching mathematical modeling for sustainability, such as Mordeson and Matthew (2020)~\cite{mordeson2021sustainable}, De Lara and Doyen (2008)~\cite{delara2008sustainable}, Gupta et al. (2016)~\cite{gupta2017mathematical}, Parkhurst (2006)~\cite{parkhurst2007introduction}, Walter (2011)~\cite{walter2011mathematics}, Fusaro and Kenschaft (2020)~\cite{fusaro2020environmental}, and Takeuchi et al. (2007)~\cite{takeuchi2007mathematics}, among others. In what follows, we will focus on reviewing Roe et al. (2018)~\cite{roe2018mathematics}.

The textbook belongs to the Texts for Quantitative Critical Thinking (TQCT) series, a collection of undergraduate textbooks that develop quantitative skills and critical thinking through the exploration of real-world questions using mathematics and statistics, providing students from various disciplines the opportunity to enhance their understanding, evaluation, and communication of quantitative information. In the foreword of the book, the former president of the Mathematical Association of America and Benediktsson-Karwa Professor of Mathematics at Harvey Mudd College, Francis Edward Su, wrote that the book presents a unique perspective on mathematics, emphasizing active participation, wisdom, and considering human factors. It explores how mathematics can contribute to addressing sustainability challenges that encompass economic, moral, and scientific aspects. As the authors have dedicated themselves to this impactful book, they hope that readers will approach it with enthusiasm to address real-world problems for the benefit of all.

This book is divided into three parts. The first part is the core of the material, where fundamental concepts are presented in the order of complexity. It comprises six chapters that cover the various topics of measuring, flowing, connecting, changing, risking, and deciding. The second part contains a collection of nine case studies that apply the mathematical tools from the previous part to address sustainability-related questions and explore real-world examples. The third part contains reference materials, including suggestions for further reading.

The authors have done an excellent job in developing the theme by building progressive themes, from measuring to deciding. However, further investigation of the mathematical topics presented in each chapter may reveal different structures. Not all users, instructors and students alike, will find the mathematical content progressive. By selecting different mathematical topics, the core of the book looks like a salad-style of subjects, a mix of everything. Instead of focusing on a particular topic, let say, differential equations, and building up from there, it requires a multidisciplinary background approach. Any instructor who plans to use it for a course should have a broad knowledge of various topics in mathematics, in addition to being aware of issues in sustainability and the ethics surrounding its implementation. A good knowledge of economics would also be advantageous, as the market paradigm, Pareto efficiency, market failure, expected utility theory, and prospect theory, among others, occur sporadically among the texts.

To provide some ideas, the second chapter on flowing requires solid knowledge of system theory, including stocks, flows, equilibrium, feedback loops, their relationships, and their dynamics. The theory also finds applications in business, management, industrial processes, organizations, and engineering, among others~\cite{sterman2000system,roberts1978managerial,bala2017system,ogata2004system,meadows2008thinking}. The third chapter on connecting requires graph theory~\cite{easley2012networks,henning2022graph,vansteen2010graph}. Both Chapters 3 and 4 on connecting and changing discussed a commonly utilized model in population growth, that is, exponential and logistic, albeit in the absence of calculus and differential equations. The fifth chapter on risking requires a profound understanding of data analysis, probability, and statistics~\cite{anderson2017introduction,bertsekas2008introduction,blitzstein2019introduction,grinstead1997introduction,morin2016probability}. Finally, the sixth chapter on deciding requires a solid understanding of game theory~\cite{brams2011game,moris1994introduction,osborne2004an,rasmusen2001games,stahl1999a,tadelis2013game}. 

It is commonly found in many textbooks on mathematical modeling that growth and decay modeling often employs initial value problems of an ordinary differential equation, a subject that requires understanding and application of calculus. Amazingly, Roe et al. (2018) presented the material in the absence of calculus. However, they allowed instructors to have any liberty in inserting calculus and differential equations when discussing this subject. Modeling population growth has direct implications for sustainability. Similar to the two sides of coins, population dynamics are closely linked to sustainability. On the one hand, an increasing population can impact on the availability of resources and the environment because there would be a greater demand for food, water, energy, and other resources that can strain ecosystems and lead to resource depletion~\cite{hunter2000the}. On the other hand, the low birth rates that are currently occurring in many parts of the world can pose challenges for sustainable development due to the declining workforce and aging population that impact social safety net, healthcare systems, labor shortages, decrease in productivity, and hindering economic growth~\cite{chan2022an,uemura2014population}.

Regarding the modeling of the decay rate, Roe et al. (2018) employed a discrete version of Newton's Law of Cooling, although their presented example does not really concern with sustainability. We would like propose the following example on this topic. A similar case in point or other real-life applications of the law associated with sustainability will enhance students' understanding of the mathematics of sustainable development. 

Farmers were faced with the task of preserving a large quantity of harvested vegetables to prevent spoilage. The vegetables had an initial temperature of 25°C, and they needed to be stored in a cold storage facility maintained at a constant temperature of 5°C. The cooling process followed Newton's Law of Cooling, with a cooling rate of 0.099°C per hour. Determine the time it would take for the vegetable temperature to reach 10°C. Additionally, if the initial temperature of the vegetables was lowered to 20°C and the cooling rate was increased to 0.183°C per hour, calculate the time difference in electricity savings compared with the previous scenario.

To solve this problem, we acquired Newton's Law of Cooling, which can be formulated as the following differential equation:
\begin{equation*}
\frac{dT}{dt} = k \left(T - T_s \right),
\end{equation*}
where $T$ denotes the temperature of the object at time $t$, $T_s$ denotes the temperature of the surroundings, and $k$ is the cooling rate constant. It solutions is given by
\begin{equation*}
T\left(t \right) = T_s + \left(T_0 - T_s \right) e^{- k t},
\end{equation*} 
where $T_0$ denotes the initial temperature of the object. 
We have $T_s = 5$, $T_0 = 25$, $k = 0.099$, and $T(t) = 10$, then
\begin{align*}
	10 &= 5 + \left(25 - 5 \right) \ e^{-0.099 \ t} \\
	5 &= 20 \ e^{-0.099 \ t} \\
e^{-0.099 \ t} &= \frac{5}{20} = \frac{1}{4} \\
	0.099 t &= \ln 4 \\
	t &= \frac{\ln 4}{0.099} \approx 14.
\end{align*}
Thus, it would take more than 14 hours for the vegetables to reach 10°C. For the second scenario, we have $T_s = 5$, $T_0 = 20$, $k = 0.183$, and $T(t) = 10$, then
\begin{align*}
	10 &= 5 + \left(20 - 5 \right) \ e^{-0.183 \ t} \\
	5 &= 15 \ e^{-0.183 \ t} \\
e^{-0.183 \ t} &= \frac{5}{15} = \frac{1}{3} \\
	0.183 t &= \ln 3 \\
	t &= \frac{\ln 3}{0.183} \approx 6.
\end{align*}
Therefore, by decreasing the initial temperature and increasing the cooling rate, it would take at least 6 hours for the vegetables to reach 10°C. Hence, the time difference in saving electricity compared with the previous case would be $14 - 6 = 8$~hours. 

This example illustrates a simple application of Newton's Law of Cooling in sustainability, that is, the food preservation and storage of perishable goods. By understanding the cooling rates and temperature changes, optimal storage conditions can be determined to extend the shelf life and reduce food waste~\cite{konovalenko2021real}. Other applications are certainly available, although they might involve and be in combination with other physical principles, such as building energy efficiency and thermal comfort analysis. Understanding how buildings lose heat in the environment is crucial for a sustainable design. Newton's Law of Cooling can be used to model and optimize heating and cooling systems, improve energy efficiency, and reduce greenhouse gas emissions~\cite{anandh2014sustainable}. The law can also help assess and ameliorate the thermal comfort of indoor spaces. By analyzing the cooling rates and temperatures in different areas, sustainable building designs can be developed to provide comfortable environments while minimizing the energy consumption~\cite{utkarsh2021ambient}.

In terms of probability, Roe et al. (2018) presented examples that are commonly found in any textbook discussing the topic, such as tossing a coin, rolling dice, drawing cards from a deck, and other similar experiments. Although one example of inferential statistical analysis is closely related to sustainability, that is, groundwater contamination by methane either due to its natural occurrence or via hydraulic fracturing~\cite{molofsky2013evaluation,osborn2011methane}, there were no particular examples on conditional probability in sustainability-related contexts. Many readers would appreciate more examples of the application of conditional probability to various scenarios related to sustainability. For example, in urban planning and transportation, conditional probability can be utilized to assess the likelihood of achieving sustainable transportation outcomes, such as reduced congestion or increased use of public transit, based on factors including infrastructure development, policy interventions, and behavior change~\cite{jiang2022sustainability,rahman2022towards}. In the area of natural resource management, conditional probability can be used in analyzing the likelihood of sustainable resource extraction, such as by calculating the probability of a fish population recovering to a desired level given certain conservation measures and fishing quotas~\cite{holland1996marine,holland1999an,punt1997fisheries}.

We consider the following example has a direct application to sustainability and simultaneously check student's understanding of conditional probability.
\begin{figure}[h]
\centering
\includegraphics[width=0.9\linewidth]{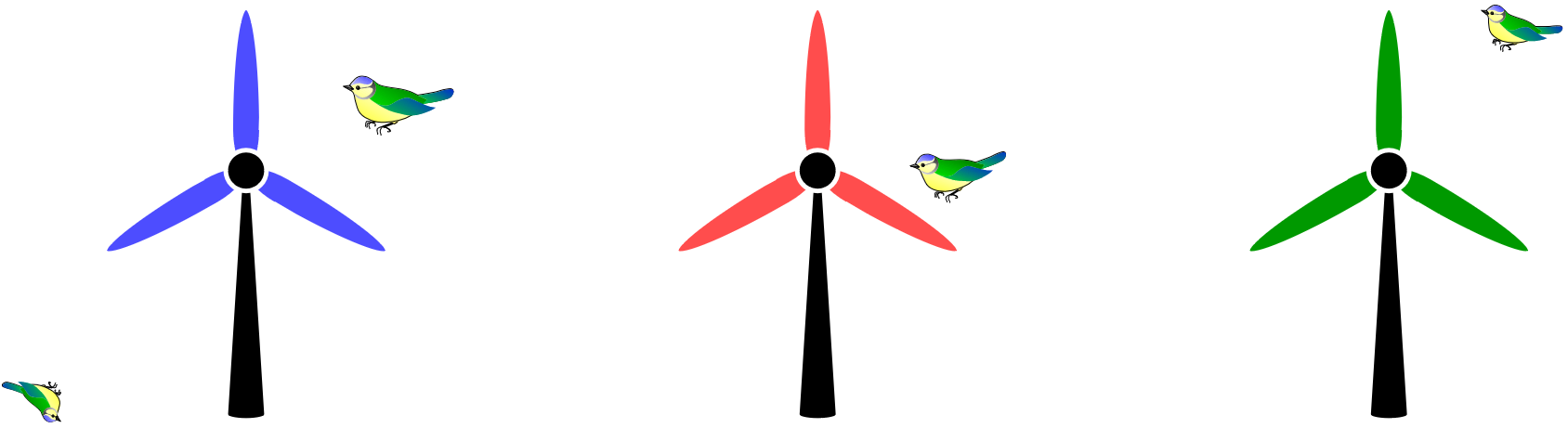}
\caption{Improper placement of wind turbines can result in unintended consequences for bird and bat fatalities, despite their significance as a valuable source of intermittent renewable energy.}  \label{fig:windturbine}
\end{figure}

The construction of a wind power infrastructure has been identified as a cause of bird and bat mortality. See Figure~\ref{fig:windturbine}. This is because turbines are often situated in the migration paths of birds, posing a risk to birds with their spinning blades. Additionally, the presence of wind turbines can attract bats owing to increased insect density, as they perceive the sites as a potential food source. However, it is worth noting that the pesticides used in agricultural areas also pose a significant threat to bird and bat populations. In a specific region examined by environmental activists, the annual mortality rates of migrating birds in the presence of wind turbines and pesticides were recorded as 27,000 and 41,000, respectively. The activists discovered that the spinning blades of the wind farm caused 9,000 bat fatalities per year, whereas agricultural pesticides accounted for 23,000 bat mortalities annually in the area. What is the probability that pesticides cause animal mortality, given that the flying creatures are birds? Find the probability of the flying creatures are bats, given that the cause of fatality is wind turbines. 

To answer these questions, constructing a two-way frequency table would definitely be helpful, as presented in Table~\ref{batbird}. Otherwise, we simply use the formula for conditional probability. Given that the flying creatures are birds, the probability that pesticides cause animal mortality is given as follows:
\begin{equation*}
P\left(\text{pesticides} \ | \ \text{bird} \right) = \frac{41,000}{41,000 + 27,000} = \frac{41}{68} = 60.29\%.
\end{equation*}
Therefore, the probability that agricultural pesticides cause animal mortality in the examined region, given that the animals are birds, is approximately 60\%. 
If the wind turbines cause animal mortality, the probability that flying creatures are bats is given as follows:
\begin{equation*}
P\left(\text{bat} \ | \ \text{wind turbines} \right) = \frac{9,000}{9,000 + 27,000} = \frac{9}{36} = 25\%.
\end{equation*}
Hence, the probability of animal mortality is bats, assuming that they were killed by wind turbines in the examined region was 25\%.
\begin{table}[h]
\begin{center}
\begin{tabular}{@{}llccr@{}}
\toprule
									&				& \multicolumn{2}{c}{Mortality causes and figures} 	& \multirow{3}{*}{Total}		 		\\ 
									&				& \multicolumn{2}{c}{(in thousands)} 				&		 		\\ \cline{3-4}									
									& 				& Wind turbines  	& Pesticides 					&				\\ \cline{2-5}
\multirow{2}{*}{Flying creatures} 	& Birds			& 27				& 41 							& 68			\\
									& Bats			& {\:}9				& 23							& 32			\\ \cline{3-5}
\multicolumn{2}{c}{Total}							& 36				& 64							& 100			\\ \bottomrule
\end{tabular}
\end{center}
\caption{Two-way frequency table for bird and bat fatalities due to wind turbines and pesticides. These figures are useful for solving problems related to the conditional probability.}		\label{batbird}
\end{table}

Although the geographical location in this example is unspecified and the figures are fictitious, this problem can serve the purpose in raising an awareness of sustainability. On the one hand, the construction of wind turbines contributes to an additional source of renewable energy, in this case, using the wind's kinetic energy to generate electricity. On the other hand, constructing such wind turbines and locating them in correct locations concern not only better energy efficiency but also environmental and social problems. Locations with the best wind currents often times coincide with migratory paths of birds~\cite{barrios2004behavioural,saidur2011environmental,welch2009the}. Additionally, research on ecological sustainability suggests that bat fatalities at wind turbines can be caused by the proximate category (direct means of death, such as collision with towers and rotating blades) and ultimate category (collisions due to bats' attraction to turbines). Possible explanations for the latter propose that bats could be drawn to turbines because of curiosity, misinterpretation, or the potential for activities such as feeding, roosting, flocking, and mating~\cite{arnett2013thresholds,cryan2009causes,frick2020a}. It turns out that bird and bat fatalities caused by wind turbines are just on a smaller scale in comparison to other causes. Communication towers, automobiles, pesticides, buildings, and cats killed more birds than wind turbines in the US alone~\cite{erickson2005a,fws,geluso1976bat,manville2009towers,mineau2002estimating,torquetti2021exposure,yang2021unprecedented}. Hence, adding pesticides to the considered example provides a better perspective for other causes of avian and bat casualties.

Regarding game theory and the tragedy of the commons, Roe et al. (2018) presented some examples related to sustainability, that is, the pollution control of a system of lakes and rivers between two territories and two farmers who share a common grassland and decide whether to add cows, respectively. The former is a bit complicated because the two countries share three lakes. A simplified example would be helpful for newcomers to game theory. The latter is an excellent one because it not only employs some concepts in game theory but also allows for some variations of different scenarios. The following example provides a simpler version of game theory on sustainability. 

Consider two small towns, Arang (A) and Bimo (B), that share a common border. They are currently evaluating two waste management options: one that focuses on minimizing an environmental impact (E) and another that aims to maximize cost-effectiveness (C). The town councils of Arang and Bimo have calculated the annual savings associated with each option, depending on the priorities that each town chooses to adopt. On the one hand, if both towns prioritize minimizing the environmental impact, Arang and Bimo would save approximately \$10,000 and \$8,000 annually, respectively. On the other hand, if both towns prioritize maximizing cost-effectiveness, Arang and Bimo would save approximately \$8,000 and \$7,000 annually, respectively. In the scenario where one town dedicates itself exclusively to either minimizing environmental impact or maximizing cost-effectiveness, each town would save approximately \$6,000 annually. However, due to the difference in size between the two towns, Arang would only save \$5,000 if Bimo town focuses on maximizing cost-effectiveness, while Arang prioritizes minimizing the environmental impact. Conversely, Bimo would save an additional \$4,000 compared with Arang if the roles were reversed. By carefully analyzing these savings and considering their different priorities, Arang and Bimo can make informed decisions on the most suitable waste management option that aligns with their sustainability goals and economic considerations. Construct a payoff matrix to represent the decision-making process for the ``game'' of selecting a waste management system between two ``players'' of towns and determine the Nash equilibrium in this waste management system game by analyzing the optimal strategies for each pair of players.

Consider the towns of Arang and Bimo as the two players in a game. Their objective was to select the most suitable waste management system for their respective towns. The payoff matrix, representing the benefits or saving costs associated with different choices, is presented in Table~\ref{arangbimo}. The given figures are in thousands of dollars. %given as follows (in thousands dollars):
\begin{table}[h]
\begin{center}
\begin{tabular}{@{}cccc@{}}
\toprule
						&				& \multicolumn{2}{c}{Bimo} 			 		\\ \cline{3-4}
						& 				& Environment 	& Cost 	\\ \cline{2-4}
\multirow{2}{*}{Arang} 	& Environment	& (10,8)		& (5,6) \\
						& Cost			& ( 6,9)		& (8,7)	\\ \bottomrule
\end{tabular}
\end{center}
\caption{The payoff matrix for the game theory problem involving two towns, Arang and Bimo, that need to find a balance between minimizing the environmental impact and maximizing cost-effectiveness.}		\label{arangbimo}
\end{table}

To determine a Nash equilibrium, we must find each player's strictly best response to the other player's strategy. We observe the following situations:
\begin{itemize}[leftmargin=1em]
\item[$\bullet$] If Arang plays the environment, then Bimo's payoffs for the environmental impact and cost-effectiveness would be \$8,000 and \$6,000, respectively. Thus, Bimo plays the environment is the strictly best response to this strategy.
\item[$\bullet$] If Arang plays cost, then Bimo's payoffs for the environment and cost are \$9,000 and \$7,000, respectively. Thus, Bimo plays environment is also still the strictly best response to this strategy.
\end{itemize}
It turns out that playing the environment is a strictly dominant strategy for Bimo. Conversely, we note the following observations:
\begin{itemize}[leftmargin=1em]
\item[$\bullet$] If Bimo plays the environment, then Arang's payoffs for the environmental impact and cost-effectiveness are \$10,000 and \$6,000, respectively. Thus, Arang plays environment is the strictly best response to this strategy.
\item[$\bullet$] If Bimo plays the cost, then Arang's payoffs for the environment and cost are \$5,000 and \$8,000, respectively. Thus, Arang plays cost turns out the strictly best response to this strategy.
\end{itemize}
We observe that Arang does not have any strictly dominant strategy in this game. However, because playing the environment is the strictly best response from each town to the other town's strategy, both Arang and Bimo play the environmental impact is a Nash equilibrium.

There are two remarks that we need to state. First, the word ``approximately,'' which appears twice in the problem, does not really mean anything, and it can be dropped as well if the readers wish. What we meant with this word is, for example, that \$9,990 is approximately \$10,000, in the same manner as \$10,050. However, figures \$4,000, \$5,000, and \$6,000 are on different orders, and they are not approximately the same values. Second, the sentence ``Conversely, Bimo would save an additional \$4,000 compared to Arang if the roles were reversed.'' might be confusing to interpret by non-native English speakers. This statement should be understood from the preceding sentences. The sentence ``In the scenario where one town dedicates itself exclusively to either minimizing environmental impact or maximizing cost-effectiveness, each town would save approximately \$6,000 annually.'' means $\text{AC} = 6$ and $\text{BC} = 6$. The next sentence ``However, due to the difference in size between the two towns, Arang would only save \$5,000 if Bimo town focuses on maximizing cost-effectiveness while Arang prioritizes minimizing environmental impact.'' means $(\text{AE}, \text{BC}) = (5,6)$. Hence, if the roles are reversed, it means $(\text{AC}, \text{BE}) = (6, 5 + 4) = (6, 9)$ and not $(\text{AC}, \text{BE}) = (6, 6 + 4) = (6, 10)$ because the latter argument might appear from reasoning that reversing the role of BC to BE would add an additional \$4,000.

This example illustrates some basic concepts in game theory: constructing a payoff matrix of a game, understanding each player's strategy and payoff relative to the other player, and finding the best response strategy and Nash equilibrium. Indeed, game theory has been applied to various aspects of sustainability, including, but not limited to, stakeholders and company leaders~\cite{lozano2011addressing}, preserving energy resources~\cite{dolinsky2015sustainable}, resource depletion strategy~\cite{seifi2018using}, decision making in the chemical industry~\cite{gonzales2021game}, and agricultural supply chain~\cite{song2022blockchain}.

The following and our final example illustrates the mathematical modeling of the tragedy of the commons~\cite{hardin1968the,hardin1974living,hardin1998extensions}. Although it lacks an application of game theory, it employs the power of an exponential model. This problem was adapted from the TED-Ed video lesson by Nicholas Amendolare~\cite{amendolare2017teded}.
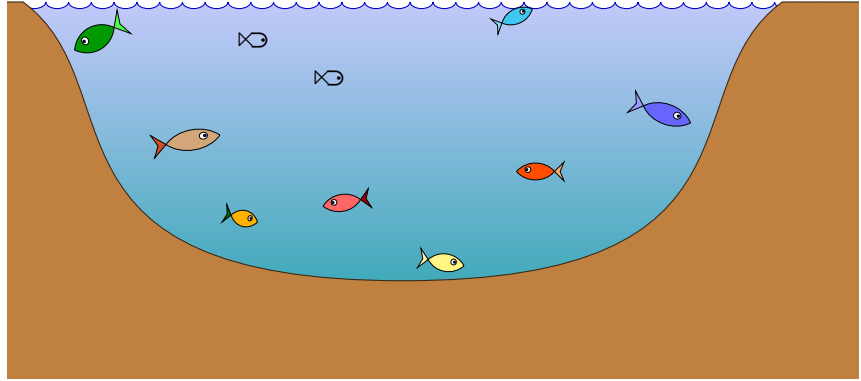
\begin{figure}[h]
\begin{center}
\begin{tikzpicture}
% Define some reference points 
% The figure is drawn a bit bigger, and then clipped to the following dimensions:
\coordinate (clipping area) at (11, 5);
\clip (-0.2,0) rectangle (clipping area);
	
% Next reference points are relative to the lower left corner of the clipping area
\coordinate (water level) at (0, 5);
\coordinate (bottom)      at (5, 1.3);     % (bottom of the pit)
\coordinate (ground1)     at (0, 5);       % (left shore)
\coordinate (ground2)     at (10, 5);      % (right shore)
	
% Coordinates of the bigger area really drawn
\coordinate (lower left)  at ([xshift=-5mm, yshift=-5mm] 0,0);
\coordinate (upper right) at ([xshift=5mm,  yshift=5mm] clipping area);
	
% Draw the water and ripples
\draw [draw=blue!80!black, decoration={bumps, mirror, segment length=6mm}, decorate, bottom color=cyan!60!black, top color=blue!20!white] 
(lower left) rectangle (water level-|upper right);
	
% Draw the ground
\draw [draw=brown!30!black, fill=brown] 
(lower left) -- (lower left|-ground1)  --
(ground1) .. controls ($(ground1)!.3!(bottom)$) and (bottom-|ground1) ..
(bottom) .. controls (bottom-|ground2) and ($(ground2)!.3!(bottom)$) .. 
(ground2) -- (ground2-|upper right) -- (lower left-|upper right) -- cycle;
	
% Draw the fish
\node[rotate=180] at (4,4) {${\cdot}\kern -4pt {\subset}\kern -3pt {\rtimes}$};
\node[rotate=180] at (3,4.5) {${\cdot}\kern -4pt {\subset}\kern -3pt {\rtimes}$};
\scalebox{0.25}{
\begin{scope}[shift={(15,8)}, rotate=10]
\draw[fill=red!60!white] (1,1) to[bend left=50] (3,1) to[bend left=50] (1,1);
\draw[fill=red!60!black] (3,1) -- (3.5,1.5) -- (3.3,1) -- (3.5,.5) -- cycle;
\draw[fill=white] (1.6,1.1) circle (.15cm); 
\draw[fill=blue] (1.55,1.1) circle (.05cm);
\end{scope}
\begin{scope}[shift={(25,10)}]
\draw[fill=orange!60!red] (1,1) to[bend left=50] (3,1) to[bend left=50] (1,1);
\draw[fill=orange!60!white] (3,1) -- (3.5,1.5) -- (3.3,1) -- (3.5,.5) -- cycle;
\draw[fill=white] (1.6,1.1) circle (.15cm); 
\draw[fill=blue] (1.55,1.1) circle (.05cm);
\end{scope}
}
\scalebox{0.3}{
\begin{scope}[shift={(10,10)}, rotate=10, xscale=-1.2]
\draw[fill=brown!70!white] (1,1) to[bend left=50] (3,1) to[bend left=50] (1,1);
\draw[fill=brown!60!red] (3,1) -- (3.5,1.5) -- (3.3,1) -- (3.5,.5) -- cycle;
\draw[fill=white] (1.6,1.1) circle (.15cm); 
\draw[fill=blue] (1.55,1.1) circle (.05cm);
\end{scope}
\begin{scope}[shift={(2,13)}, rotate=30, yscale=1.2]
\draw[fill=green!60!black] (1,1) to[bend left=50] (3,1) to[bend left=50] (1,1);
\draw[fill=green!70!white] (3,1) -- (3.5,1.5) -- (3.3,1) -- (3.5,.5) -- cycle;
\draw[fill=white] (1.6,1.1) circle (.15cm); 
\draw[fill=blue] (1.55,1.1) circle (.05cm);
\end{scope}
\begin{scope}[shift={(30,10)}, rotate=-20, xscale=-1.1]
\draw[fill=blue!60!white] (1,1) to[bend left=50] (3,1) to[bend left=50] (1,1);
\draw[fill=blue!40!white] (3,1) -- (3.5,1.5) -- (3.3,1) -- (3.5,.5) -- cycle;
\draw[fill=white] (1.6,1.1) circle (.15cm); 
\draw[fill=blue] (1.55,1.1) circle (.05cm);
\end{scope}
}
\scalebox{0.2}{
\begin{scope}[shift={(30,6)}, rotate=-10, xscale=-1.2,yscale=1.3]
\draw[fill=yellow!60!white] (1,1) to[bend left=50] (3,1) to[bend left=50] (1,1);
\draw[fill=yellow!40!white] (3,1) -- (3.5,1.5) -- (3.3,1) -- (3.5,.5) -- cycle;
\draw[fill=white] (1.6,1.1) circle (.15cm); 
\draw[fill=blue] (1.55,1.1) circle (.05cm);
\end{scope}	
\begin{scope}[shift={(35,24)}, rotate=25, xscale=-1.1,yscale=1.1]
\draw[fill=cyan!60!white] (1,1) to[bend left=50] (3,1) to[bend left=50] (1,1);
\draw[fill=cyan!40!white] (3,1) -- (3.5,1.5) -- (3.3,1) -- (3.5,.5) -- cycle;
\draw[fill=white] (1.6,1.1) circle (.15cm); 
\draw[fill=blue] (1.55,1.1) circle (.05cm);
\end{scope}
\begin{scope}[shift={(16,9)}, rotate=-15, xscale=-0.9,yscale=1.2]
\draw[fill=orange!60!yellow] (1,1) to[bend left=50] (3,1) to[bend left=50] (1,1);
\draw[fill=green!40!black] (3,1) -- (3.5,1.5) -- (3.3,1) -- (3.5,.5) -- cycle;
\draw[fill=white] (1.6,1.1) circle (.15cm); 
\draw[fill=blue] (1.55,1.1) circle (.05cm);
\end{scope}
}	
\end{tikzpicture}
\end{center}
\caption{The illustration depicts a pond with the potential for exponential growth of the fish population. However, irresponsible fishing practices can lead to the tragedy of the commons, resulting in the depletion of fish stock.}	\label{fishpond}
\end{figure}

In Green Village, there is a pond that serves as a crucial source of livelihood for villagers, providing them with fish. See Figure~\ref{fishpond}. The number of fish in the pond doubled every month, ensuring a potential increase in catch. Alice and Bob, two residents of the village, have the responsibility of catching fish to sustain their community. They understood the importance of preserving the pond's resources and decided not to catch all the fish at once, allowing for reproduction. In the first month, both Alice and Bob catch one fish each, and plan to double their catch every subsequent month. However, despite their best intentions, they lack expertise in mathematical modeling, particularly in the field of sustainability. As a consequence of their miscalculations, the tragic phenomenon known as the ``tragedy of the commons'' became a concern when the fish stock in the pond faced depletion. If initially there are 5 pairs of fish in the pond, how many months will it take for all fish to be gone? To prevent the tragedy of the commons by the end of the year, what should be the minimum total number of fish in the pond at the beginning of the year?

To answer the first question, let $m_0 = 2q = 2^{0 + 1} (q - 0)$ be the total number of fish at the beginning of the month, where $q \in \mathbb{N}$. Let $p_n = 2^n$ be the total number of fish that were caught by Alice and Bob at the end of the month $n \in \mathbb{N}$. The following recursive relation indicates the total number of fish at the end of each month:
\begin{align*}
	m_1 &= 2 \left(m_0 - p_1 \right) = 2 \left(2q - 2 \right) = 4 \left(q - 1 \right) = 2^{1 + 1} \left(q - 1\right) \\
	m_2 &= 2 \left(m_1 - p_2 \right) = 2 \left[ 4\left(q - 1 \right) - 4 \right]  = 8 \left(q - 2 \right) = 2^{2 + 1} \left(q - 2 \right) \\
	m_3 &= 2 \left(m_2 - p_3 \right) = 2 \left[ 8\left(q - 2 \right) - 8 \right]  = 16 \left(q - 3 \right) = 2^{3 + 1} \left(q - 3 \right) \\
	& \qquad \vdots \\
	m_n &= 2 \left(m_{n - 1} - p_n \right) = 2 \left[ 2^n \left(q - n \right) - 2^n \right] = 2^{n + 1} \left(q - n \right).
\end{align*} 
Because $q = 5$ and we want to solve $m_n = 0$, then $q - n = 5 - n = 0$, which gives $n = 5$. Hence, after 5 months, all fish in the pond will disappear. Table~\ref{fishstock} shows the situation describing fish stock depletion, where the number of fish caught follows monthly exponential growth. Observe that despite the potential for exponential reproduction, the initial stock of the 5 pairs of fish would be depleted after 5 months.

To answer the second question, we use the last equation for $n = 12$. Because we want $m_{12} > 0$ or $q - 12 > 0$, so $q > 12$. Hence, considering $q = 13$, then the total number of fish at the beginning of the year should be at least 26 to avoid the tragedy of the commons after one year.
\begin{table}[h]
\begin{center}
\begin{tabular}{@{}ccc@{}}
\toprule
Month $n$ 	& Fish caught $p_n$ 	& Fish in month $n$, $m_n$ 	\\ \midrule
	0		&		0				&		10					\\
	1		&		2				&		16					\\
	2		&		4				&		24					\\
	3		&		8				&		32					\\
	4		&		16				&		32					\\	
	5		&		32				&		 0					\\	
\bottomrule	
\end{tabular}
\end{center}
\caption{Situation describing fish stock depletion. The total number of fish at the end of each month with exponential growth in catch is listed in the third column. Despite the potential for exponential reproduction, the initial stock of 10 fish would be exhausted after five months.}		\label{fishstock} %If there were 10 fish at the beginning of the month, the stock of fish would be depleted after five months even though they could reproduce double figures every month.}
\end{table}

While this simplified illustration portrays the tragedy of the commons, its occurrence is evident in the fishing industry, where fisheries frequently experience over-exploitation. When a fishery is regarded as a common resource, both locally and internationally, fishermen often prioritize maximizing their catch to serve their immediate self-interest, often disregarding long-term implications. Consequently, their unsustainable practices endanger the future sustainability of fisheries, posing a threat to the provision of a reliable food source for present and future generations. In addition, these actions contribute to the degradation of fish habitats and ecosystems~\cite{alfattal2009the,berkes1985fishermen,kraak2011exploring,leal1998community,mcwhinnie2009the}.

% Conclusion
\section{Conclusion}			\label{discussconclude}

We have conducted a comprehensive review of the current state of teaching and learning mathematical modeling at the tertiary level, with a particular emphasis on sustainability and sustainable learning in education. This resonates with UN Sustainable Development Goal (SDG) 4, which promotes continuous and lifelong learning for all individuals. Our analysis revealed a significant gap in the body of published literature regarding mathematical modeling for sustainability and sustainable learning in mathematics education.

Through our examination of various textbooks on mathematical modeling, we discovered that most of them cover a wide range of topics but lack a specific focus on sustainability. Although some textbooks include examples related to sustainability, they are often buried within the content and require additional effort from instructors to extract and tailor them to sustainability issues. However, many textbooks exhibit common characteristics, including considerations such as the type of variables (discrete or continuous), the incorporation of randomness and uncertainty in modeling (deterministic or stochastic), and the approach employed to solve problems (analytical or numerical).

To address this gap, we proposed a liberal arts course on mathematical modeling for sustainability that is suitable for students from diverse academic backgrounds. This course not only introduces students to sustainable development and the mathematics of sustainability but also enhances their quantitative thinking skills for sustainable learning in mathematics education. By bridging these two domains, the course aimed to promote a quantitative understanding of sustainability and the application of mathematical skills in various contexts. Our aspiration is for students who have completed this course to embrace sustainable learning in mathematics education, which entails a commitment to ongoing learning, relearning, and continuous improvement in mathematical skills beyond the confines of formal education.

We also explored the potential benefits of adopting a specific textbook in the course on mathematical modeling for sustainability. Although it may not cover traditional college mathematics topics such as the use of calculus and differential equations, it offers quantitative content and encourages critical thinking through real-world problem-solving. The book's accessible language appeals to students outside STEM fields, while those with a strong mathematical background can gain a broader perspective on the mathematics of sustainability. Additionally, we provided classroom examples that demonstrate how mathematical concepts, such as probability and game theory, find practical applications in sustainability.

In conclusion, our research highlights the need for a focused exploration of mathematical modeling for sustainability in higher education. By integrating sustainability into mathematics education, we can equip students with the skills and knowledge necessary for addressing complex real-world challenges in a sustainable and quantitative manner.

\subsection*{Acknowledgments}
The author acknowledges the anonymous referees for reviewing this manuscript.

\subsection*{Funding statement}
This research receives no funding.

\subsection*{Conflicts of interest}
The author declares that he has no conflicts of interest to disclose.

\subsection*{Ethics statement}
Not applicable.

% References
\end{document}